\documentclass[11pt]{article}
\usepackage{amssymb}
\usepackage{amsmath}
\usepackage{amsthm}
\usepackage{latexsym}
\usepackage{amsfonts}
\usepackage{graphicx}
\usepackage{graphics}

\newcommand{\dis}{\displaystyle}
\theoremstyle{plain}
\newtheorem{thm}{Theorem}[section]   

\newtheorem{lem}[thm]{Lemma}

\theoremstyle{definition}
\newtheorem{rem}[thm]{Remark}

\newtheorem*{Proof}{Proof}

\newcommand{\fa}{\forall}

\newcommand{\ra}{\;\rightarrow\;}

\newcommand{\de}{\delta }
\newcommand{\OO} {{\varOmega}}
\newcommand{\De} {{\varDelta}}
\newcommand{\e}{\varepsilon }

\newcommand{\la}{\lambda }

\newcommand{\ti}{\tau }

\newcommand{\C}{\mathbb{C}}

\newcommand{\N}{\mathbb{N}}

\newcommand{\ssum}{\sum\limits}

\newcommand{\ld}{\ldots}

\newcommand{\sm}{\smallsetminus}

\newcommand{\qb}{$\quad\blacksquare$}
\begin{document}
\title{\bf Generic Approximation of functions\\ by their Pad\'{e} approximants, I}
\author{G. Fournodavlos}
\date{}
\maketitle
%
\noindent
%
\begin{abstract}

Approximation of entire functions by their Pad\'{e} approximants
has been examined in the past. It is true that generically such an
approximation holds. However, examining this problem from another
viewpoint, we obtain stronger generic results on functions defined
on simply connected domains or even open sets of arbitrary
connectivity.
\end{abstract}
{\em Subject Classification MSC2010}\,: primary 41A21, 30K05
secondary 30B10, 30E10, 30K99, 41A10, 41A20. \vspace*{0.2cm} 
\section{Introduction}\label{sec1} 
\noindent

Every holomorphic function on a disc can be approximated by the
partial sums of its Taylor development. Further, in any simply
connected domain quasi all holomorphic functions are limits of a
subsequence of their partial sums in the topology of uniform
convergence on compacta (\cite{6}, \cite{7}).

Instead of considering approximation by the partial sums of the
Taylor developments, which are polynomials, one can examine the same
question using rational functions, namely the Pad\'{e} approximants
$[p/q]_f$, $p,q\ge0$ (see \cite{1}). In \cite{2} it is proved that
quasi all entire functions $f\in H(\C)$ are the limit of a
subsequence $[p_n/q_n]_f$, where $p_n-q_n\ra\infty$. Inspired by
\cite{3} we examine the same problem from a different scope and we
obtain the same result under the weaker assumption $p_n\ra+\infty$.

The same proof works if $\C$ is replaced by a simply connected open
set $\OO\subset\C$ containing 0. This is done in Section 3 of the
present article. We also mention that the condition
$p_n-q_n\ra\infty$ (or $p_n\ra+\infty$, $q_n\ra+\infty$) is used in
Theorem 5.1 of \cite{3}, and our result improves a corollary of Th.
5.1 of \cite{3}.

Finally, we mention that when we do approximation by polynomials,
the maximum principle leads us to consider compact sets $K$ with
connected complement. However, when we do approximation by rational
functions, as the Pad\'{e} approximants, we may have poles on the
holes of the compact set $K$. Thus, the result of Section 3 can be
generalized to the case of open sets $\OO\subset\C$ containing 0 of
arbitrary connectivity, under the assumption $p_n\ra+\infty$,
$q_n\ra+\infty$. This is done in Section \ref{sec4} of the present
article.

Our method of proof uses Baire's Category Theorem. We refer to
\cite{5} and \cite{4} for the role of Baire's Theorem in Analysis.
\section{Preliminaries}\label{sec2}
\noindent

Let $\OO\subseteq\C$ be an open set. Define
$K_n=\{z\in\C:\text{dist}(z,\OO^c)\ge1/n$, $|z|\le n\}$, $n\in\N$.
\begin{rem}\label{rem2.1}
i) $K_n$ is a compact subset of $\OO$ and $K_n\subseteq
K^0_{n+1}$, $\fa\;n$

ii) $\bigcup\limits^\infty_{n=1}K_n=\OO$ and if $K\subseteq\OO$
compact $\exists\;n_0\in\N$: $K_{n_0}\supseteq K$.

iii) Every component of $\C\cup\{\infty\}\sm K_n$ contains at least
one component of $(\C\cup\{\infty\})\sm\OO$ (\cite{9}).
\end{rem}

We define a metric $\rho$ on the set $H(\OO)$ (of holomorphic in
$\OO$):
\[
\rho(f,g)=\sum^\infty_{n=1}\frac{1}{2^n}\min\{\|f-g\|_{K_n},1\},
\]
where $\|\cdot\|_{K_n}$ denotes the supremum norm on $K_n$. It is
easy to see that a sequence in $H(\OO)$, $(f_m)_{m\in\N}$ converges
$f_m\overset{\rho}{\longrightarrow}f$, if and only if $f_m\ra f$
uniformly on the compact subsets of $\OO$. The space $(H(\OO),\rho)$
is a compete metric space.

Let $f$ be a function holomorphic in a neighborhood of 0 and let
$f(z)=\ssum^\infty_{v=0}a_vz^v$ its Taylor series. A Pad\'{e}
approximant $[p/q]_f$ of $f$, $p,q\in\{0,1,2,\ld\}$, is a rational
function of the form
\[
\frac{\ssum^p_{v=0}n_vz^v}{\ssum^q_{v=0}d_vz^v}, \ \ d_0=1.
\]
such that its Taylor series $\ssum^\infty_{v=0}b_vz^v$ coincides
with $\ssum^\infty_{v=0}a_vz^v$ up to the first $p+q+1$ terms;
that is $b_v=a_v$ for $v=0,\ld,p+q$ (\cite{1}).

In case of $q=0$ there exists always a unique
Pad\'{e} approximant of $f$ and $[p/q]_f(z)= S_p(z)$, where $S_p(z)= \dis\sum^p_{v=0}a_vz^v$.
For $q\ge 1$ it is true that there exists a unique Pad\'{e} approximant of $f$,
if and only if the following
determinant is not zero:
\[
\det\left|\begin{array}{cccc}
  a_{p-q+1} & a_{p-q+2} & \cdots & a_p \\
  a_{p-q+2} & a_{p-q+3} & \cdots & a_{p+1} \\
  \vdots & \vdots &  & \vdots \\
  a_p & a_{p+1} & \cdots & a_{p+q-1} \\
\end{array}\right|\neq0, \ \ a_i=0, \ \ \text{when} \ \ i<0.
\eqno(\ast)
\]
Then we write $f\in D_{p,q}$.

If $f\in D_{p,q}$, then $[p/q]_f$ $(q\ge1)$ is given by the Jacobi
explicit formula:
\[
[p/q]_f=\frac{\det\left|\begin{array}{cccc}
  z^qS_{p-q}(z) & z^{q-1}S_{p-q+1} & \cdots & S_p(z) \\
  a_{p-q+1} & a_{p-q+2} & \cdots & a_{p+1} \\
  \vdots & \vdots &  & \vdots \\
  a_p & a_{p+1} & \cdots & a_{p+q} \\
\end{array}\right|}{\det\left|\begin{array}{cccc}
  z^q & z^{q-1} & \cdots & 1 \\
  a_{p-q+1} & a_{p-q+2} & \cdots & a_{p+1} \\
  \vdots & \vdots &  & \vdots \\
  a_p & a_{p+1} & \cdots & a_{p+q} \\
\end{array}\right|},
\]
with $S_k(z)=\left\{\begin{array}{cc}
  \ssum^k_{v=0}a_vz^v, & k\ge0 \\
  0, & k<0. \\
\end{array}\right.$

If $K$ is any set we write $f\in H(K)$ if $f$ is holomorphic in
some open set containing $K$.
\begin{lem}\label{lem2.2}
Let $\OO\subseteq\C$ be an open set, $0\in\OO$, and $\la\in\N$ such
that $\overline{\De(0,r)}\subseteq K^0_\la$ for a certain $r>0$. If
$f\in H(K_\la)$, $f\in D_{p,q}$ and its Pad\'{e} approximant has no
poles in $K_\la$ and if $\e>0$ is given, then there exists $\de>0$
such that for every $g\in H(K_\la)$ with $\|g-f\|_{K_\la}<\de$ it
holds $g\in D_{p,q}$ and $\|[p/q]_g-[p/q]_f\|_{K_\la}<\e$.
\end{lem}
\begin{Proof}
Observe that the above determinant $(\ast)$ and the coefficients of
the numerator and the denominator of $[p/q]_f$ depend polynomially
on $\frac{f^{(v)}(0)}{v!}$, $v=0,1,\ld,p+q$. This implies that there
exists a $\widetilde{\de}>0$ such that for every $g\in H(K_\la)$
with
$\big|\frac{g^{(v)}(0)}{v!}-\frac{f^{(v)}(0)}{v!}\big|<\widetilde{\de}$,
$v=0,1,\ld,p+q$ it holds $g\in D_{p,q}$ and
$\|[p/q]_g-[p/q]_f\|_{K_\la}<\e$.

If $0<\de<\min\{r^v\cdot\widetilde{\de}\;|\;v=0,1,\ld,p+q\}$ and
$\|g-f\|_{K_\la}<\de$, then by Cauchy's estimates we obtain:
\[
\bigg|\frac{g^{(v)}(0)}{v!}-\frac{f^{(v)}(0)}{v!}\bigg|=
\bigg|\frac{(g-f)^{(v)}(0)}{v!}\bigg|\le\frac{\|g-f\|_{\overline{\De(0,r)}}}{r^v}\le
\frac{\|g-f\|_{K_\la}}{r^v}<\frac{\de}{r^v}<\widetilde{\de}.
\qquad \text{\qb}
\]
\end{Proof}
\begin{rem}\label{rem2.3}
It follows that $D_{p,q}$ is open in $H(\OO)$.
\end{rem}
\begin{rem}\label{rem2.4}
If all of the coefficients $\frac{f^{(v)}(0)}{v!}=a_v$,
$v=0,1,\ld,p+q$, involved in the determinant $(\ast)$ depend
linearly on $d\in\C$, $a_v=c_v\cdot d+\ti_v$, such that $c_v=0$,
when $v<p$ and $c_p\neq0$, then the determinant $(\ast)$ is a
polynomial in $d$ of order $q$ and hence only for finite values of
$d$ the determinant is zero.
\end{rem}
\section{The simply connected case}\label{sec3}
\noindent

Let $\OO\subseteq\C$ be a simply connected domain containing 0.
Also, let $F\subseteq\N\times\N$ which contains a sequence
$(\widetilde{p}_m,\widetilde{q}_m)_{m\in\N}$, such that
$\widetilde{p}_m\ra+\infty$. We define
\begin{itemize}
\item $B_F=\{f\in H(\OO)$: there exists $(p_m,q_m)_{m\in\N}$ in
$F$ such that $f\in D_{p_m,q_m}$, for all $m\in\N$ and for every
$K\subseteq\OO$ compact $[p_m/q_m]_f\ra f$ uniformly on $K\}$.
\item $E(n,s,(p,q))=\{f\in H(\OO):f\in D_{p,q}$ and
$\|[p/q]_f-f\|_{K_n}<1/s\}$, $n,s\in\N$.
\end{itemize}
\begin{lem}\label{lem3.1}
$B_F=\bigcap\limits^\infty_{n,s=1}\bigcup\limits_{(p,q)\in
F}E(n,s,(p,q)$.
\end{lem}
\begin{Proof}
It is standard and is omitted. [A similar proof can be found in
\cite{8}].
\qb
\end{Proof}
\begin{lem}\label{lem3.2}
$E(n,s,(p,q))$ is open.
\end{lem}
\begin{Proof}
$D_{p,q}$ is open (Remark \ref{rem2.3}) and similarly as in Lemma
\ref{lem2.2}, we can prove that the map
$f\mapsto\|[p/q]_f-f\|_{K_n}$ is continuous, according to the Jacobi
formula combined with Cauchy estimates. The lemma follows easily.
\qb
\end{Proof}
\begin{thm}\label{thm3.3}
$B_F$ is $G_\de$ and dense. (Hence $B_F\neq\emptyset$).
\end{thm}
\begin{Proof}
Lemma \ref{lem3.2} implies that $\bigcup\limits_{(p,q)\in
F}E(n,s(p,q))$ is open. By Lemma \ref{lem3.1} $B_F$ is $G_\de$. We
claim that $\bigcup\limits_{(p,q)\in F}E(n,s,(p,q))$ is dense. If
that is true, then Baire's Category theorem completes the proof.
By Runge's theorem the polynomials are dense in $H(\OO)$.
Therefore, it suffices to prove that for every polynomial $P$ and
$\e>0$ there exists $f\in\bigcup\limits_{(p,q)\in F}E(n,s,(p,q))$
such that $\|f-P\|_{K_N}<\e$, where $N=N(\e)\in\N$.
\begin{itemize}
\item Let $P$ be a polynomial and $(p,q)\in F$ such that
$p>\text{deg}P$.
\end{itemize}

If $q=0$, define $f(z)=P(z)+dz^p$, $d\in\C\sm\{0\}$. It is
immediate that $f\in D_{p,q}$ and $[p/q]_f=f$. It follows $f\in
E(n,s,(p,q))$. Furthermore,
$\|f-P\|_{K_N}=|d|\cdot\|z\|^p_{K_N}<\e$, when
$0<|d|<\e/\|z\|^p_{K_N}$.

If $q\ge1$, we define
$\widetilde{f}(z)=\frac{P(z)+dz^p}{1-(cz)^q}$, $d,c\in\C\sm\{0\}$,
where $c$ and $d$ will be determined later on.
\begin{itemize}
\item Let $\la\in\N$, $\la>\max\{n,N\}$, such that
$K^0_\la\supseteq\overline{\De(0,r)}$, where $r>0$ is fixed.
\item We have
$\dis\inf_{K_\la}|1-(cz)^q|\ge1-|c|^q\cdot\|z\|^q_{K_\la}>\frac{1}{2}$,
when $0<|c|<\big(\frac{1}{2\cdot\|z\|^q_{K_\la}}\big)^{1/q}$.
\item There exists $\widetilde{\de}>0$ such that
$\|\widetilde{f}-P\|_{K_\la}=\dis\sup_{z\in
K_\la}\frac{|(cz)^q\cdot P(z)+d\cdot z^p|}{|1-(cz)^q|}\le
2(|c|^q\cdot\|z\|^q_{K_\la}\cdot\|P(z)\|_{K_\la}+|d|\cdot\|z\|^p_{K_\la})
<\e/2$, whenever $|d|<\widetilde{\de}$ and
$|c|<\widetilde{\de}<\big(\frac{1}{2\|z\|^q_{K_\la}}\big)^{1/q}$.
\item Around 0,
$\widetilde{f}(z)=P(z)+dz^p+P(z)\cdot(cz)^q+dz^p\cdot(cz)^q+\cdots$.
We fix a constant $c$ satisfying the above. According to Remark
\ref{rem2.4} we can choose $0<|d|<\widetilde{\de}$ such that
$\widetilde{f}\in D_{p,q}$. By the uniqueness of the Pad\'{e}
approximant of $\widetilde{f}$ we obtain
$[p/q]_{\widetilde{f}}=\widetilde{f}$.
\item Applying Lemma \ref{lem2.2} there exists
$0<\de<\min\{1/2s,\e/2\}$ such that, if $f\in H(K_\la)$ with
$\|\widetilde{f}-f\|_{K_\la}<\de$ it follows $f\in D_{p,q}$ and
$\|[p/q]_{\widetilde{f}}-[p/q]_f\|_{K_\la}<1/2s$. By Runge's
theorem we can choose $f$ to be a rational function with poles
only in $(\C\cup\{\infty\})\sm K_\la$. More particularly Remark
\ref{rem2.1} enables us to choose $f$ with pole only at $\infty$
because $(\C\cup\{\infty\}\sm\OO$ is connected. Thus, $f$ is a
polynomial and $f\in H(\OO)$. We also have
\item
$\|[p/q]_f-f\|_{K_n}\le\|[p/q]_f-f\|_{K_\la}\le\|[p/q]_f-[p/q]_{\widetilde{f}}\|_{K_\la}
+\|\widetilde{f}-f\|_{K_\la}<1/2s+\de<1/s$. It follows that
$f\in E(n,s,(p,q))$ and
$\|f-P\|_{K_N}\le\|f-P\|_{K_\la}\le\|f-\widetilde{f}\|_{K_\la}+\|\widetilde{f}-P\|_{K_\la}<
\de+\frac{\e}{2}<\e$.
\end{itemize}

This completes the proof.  \qb
\end{Proof}
\begin{rem}\label{rem3.4}
In the above proof we have not used the fact that $\OO$ is
connected. Therefore Theorem \ref{thm3.3} is valid for any simply
connected open set $\OO\subseteq\C$ containing 0.
\end{rem}
\section{Domains of arbitrary connectivity}\label{sec4}
\noindent

Let $\OO\subseteq\C$ be an open set containing 0. Also, let
$F\subseteq\N\times\N$ containing a sequence
$(\widetilde{p}_m,\widetilde{q}_m)_{m\in\N}$ such that
$\widetilde{p}_m\ra+\infty$ and $\widetilde{q}_m\ra+\infty$. We
define $B_F$ and $E(n,s,(p,q))$ similarly as in Section
\ref{sec3}.

The analogue of Lemmas \ref{lem3.1}, \ref{lem3.2} hold in this
case also and the proofs are similar.
\begin{thm}\label{thm4.1}
$B_F$ is $G_\de$ and dense. (Hence $B_F\neq\emptyset$).
\end{thm}
\begin{Proof}
Since $\bigcup\limits_{(p,q)\in F}E(n,s,(p,q))$ is open, it follows
that $B_F$ is $G_\de$. In order to use Baire's Category theorem we
will prove that $\bigcup\limits_{(p,q)\in F}E(n,s,(p,q))$ is dense.
By Runge's theorem the rational functions with poles in
$\OO^c\cup\{\infty\}$ are dense in $H(\OO)$. Therefore, it suffices
to prove that for every rational function $R$ with poles off $\OO$
and $\e>0$ there exists an $f\in\bigcup\limits_{(p,q)\in
F}E(n,s,(p,q))$ such that $\|f-R\|_{K_N}<\e$, where $N=N(\e)\in\N$.
\begin{itemize}
\item Let $R(z)=\frac{A(z)}{B(z)}$ be a rational function with
poles only in $(\C\cup\{\infty\})\sm\OO$ where $A$, $B$ are
polynomials. There is $(p,q)\in F$ such that $p>\text{deg}A$ and
$q>\text{deg}B$. We define
$\widetilde{f}(z)=\frac{A(z)+dz^p}{B(z)-(cz)^q}$,
$d,c\in\C\sm\{0\}$ where $c$ and $d$ will be determined later on.
Let $\la\in\N$, $\la>\max\{n,N\}$ such that
$K^0_\la\supseteq\overline{\De(0,r)}$ for some $r>0$.
\item Since $R$ has no poles in $\OO$, it follows $B(0)\neq0$ and
$\dis\inf_{z\in K_\la}|B(z)|>0$.
\item We have $\dis\inf_{z\in K_\la}|B(z)-(cz)^q|\ge\dis\inf_{z\in
K_\la}|B(z)|-|c|^q\cdot\|z\|^q_{K_\la}>0$, when
$0<|c|<\frac{(\dis\inf_{z\in K_\la}|B(z)|)^{1/q}}{\|z\|_{K_\la}}$.
\item There exists $\widetilde{\de}>0$ such that
$\|\widetilde{f}-R\|_{K_\la}\le\frac{\|A(z)\|_{K_\la}\cdot|c|^q\cdot\|z\|^q_{K_\la}+|d|\cdot
\|z\|^p_{K_\la}\cdot\|B(z)\|_{K_\la}}{\dis\inf_{z\in
K_\la}|B(z)|\cdot\dis\inf_{z\in K_\la}|B(z)-(cz)^q|}<\e/2$
whenever $|d|<\widetilde{\de}$ and
$|c|<\widetilde{\de}<\frac{(\dis\inf_{z\in
K_\la}|B(z)|)^{1/q}}{\|z\|_{K_\la}}$.
\item Around 0 we have: $\widetilde{f}(z)=B^{-1}(0)\cdot
A(z)+B^{-1}(0)\cdot
dz^p-B^{-1}(0)A(z)\cdot(B^{-1}(0)\widetilde{B}(z)-1)-dz^pB^{_1}(0)(B^{-1}(0)\widetilde{B}(z)-1)+
\cdots$ where $\widetilde{B}(z)=B(z)-(cz)^q$. We fix $c$
satisfying the above. Then by Remark \ref{rem2.4} we can choose
$0<|d|<\widetilde{\de}$ such that $\widetilde{f}\in D_{p,q}$.
Thus, there exists a unique Pad\'{e} approximant of
$\widetilde{f}$ and
\[
\widetilde{f}(z)=\frac{B^{-1}(0)\cdot
A(z)-B^{-1}(0)dz^p}{B^{-1}(0)\cdot B(z)-B^{-1}(0)(cz)^q} \ \
\text{satisfies} \ \ [p/q]_{\widetilde{f}}=\widetilde{f}.
\]
\item Lemma \ref{lem2.2} provides $0<\de<\min\{1/2s,\e/2\}$ such
that for every $f\in H(K_\la)$ with
$\|\widetilde{f}-f\|_{K_\la}<\de$ it follows $f\in D_{p,q}$ and
$\|[p/q]_{\widetilde{f}}-[p/q]_f\|_{K_\la}<1/2s$. By Runge's
theorem there exists a rational function $f$ with poles only in
$(\C\cup\{\infty\})\sm K_\la$ which satisfies
$\|\widetilde{f}-f\|_{K_\la}<\de$. The Remark \ref{rem2.1} allows
us to choose the poles in $(\C\cup\{\infty\})\sm\OO$ and hence
$f\in H(\OO)$.
\item We have
$\|[p/q]_f-f\|_{K_n}\le\|[p/q]_f-f\|_{K_\la}\le\|[p/q]_f-[p/q]_{\widetilde{f}}
\|_{K_\la}+\|\widetilde{f}-f\|_{K_\la}<1/2s+\de<1/s$. It follows
$f\in E(n,s,(p,q))$ and
$\|f-R\|_{K_N}\le\|f-R\|_{K_\la}\le\|f-\widetilde{f}\|_{K_\la}+\|\widetilde{f}-R\|_{K_\la}<
\de+\e/2<\e$.  \qb
\end{itemize}
\end{Proof}
{\bf Acknowledgement}: Completing this work, I would like to
express my warmest thanks and gratitude to Professor V. Nestoridis
for his guidance and his valuable instruction.
\bigskip
Department of Mathematics, \\
University of Athens \\
Panepistemiopolis, 157 84 Athens, Greece \\
e-mail: gregdavlos@hotmail.com


\begin{thebibliography}{99}
\bibitem{1} G. A. Baker, Jr. and P. R. Graves-Morris: Pad\'{e}
Approximants, Vol. 1 and 2, (Encyclopedia of Mathematics and its
Applications), Cambridge Un. Press, 2010.
%
\bibitem{2} P. B. Borwein: The usual behaviour of rational
approximants, Canad. Math. Bull. Vol. 26, 1983, p. 317-323.
%
\bibitem{3} N. J. Daras and V. Nestoridis: Universal Pad\'{e}
approximation, arXiv: 1102.4782v1, 2011.
%
\bibitem{4} K.-G. Grosse-Erdmann: Universal families and
hypercyclic operators, Bull. Amer. math. Soc., Vol. 36, 1999, p.
345-381.
%
\bibitem{5} J.-P. Kahane: Baire's Category theorem and
trigonometric series, J. Anal. Math., Vol. 80, 2000, p. 143-182.
%
\bibitem{6} W. Luh: Universal approximation properties of
overconvergent power series on open sets, Analysis, Vol. 6,
1986, p. 191-207.
%
\bibitem{7} A. Melas and V. Nestoridis: Universality of Taylor
series as a generic property of holomorphic functions, Adv. Math.
Vol. 157, 2001, p. 138-176.
%
\bibitem{8} V. Nestoridis: Universal Taylor series, Annales de l'
Institute Fourier, Vol. 46, 1996, p. 1293-1306.
%
\bibitem{9} W. Rudin, Real and Complex Analysis, McGraw-Hill,
1986.
\end{thebibliography}
\end{document}